\documentclass[pdflatex,sn-mathphys]{sn-jnl}

\usepackage{graphicx}
\usepackage{latexsym}
\usepackage{enumerate}
\usepackage{amsmath}
\usepackage{amsfonts}
\usepackage{amssymb}
\usepackage{amsthm}
\usepackage{amsopn}
\usepackage{amstext}
\usepackage{amscd}

\theoremstyle{thmstyleone}%
\theoremstyle{thmstyletwo}%

\theoremstyle{thmstylethree}%

\newcommand{\be}{\begin{equation}} 
\newcommand{\ee}{\end{equation}}
\newcommand{\beq}{\begin{eqnarray}}
\newcommand{\eeq}{\end{eqnarray}}
\newcommand{\nbeq}{\begin{eqnarray*}}
\newcommand{\neeq}{\end{eqnarray*}}

\newcommand{\alphab}{{\pmb \alpha }}

\begin{document}

\title[On characterization of the exponential distribution]
{On characterization of the exponential distribution via hypoexponential distributions  }

\author*[1,2]{\fnm{George} \sur{Yanev}}\email{george.yanev@utrgv.edu}

\affil*[1]{\orgdiv{School of Mathematical and Statistical Sciences}, \orgname{University of Texas Rio Grande Valley}, \orgaddress{\street{1201 W. University Dr.}, \city{Edinburg} \postcode{78539}, \state{Texas}, \country{USA}}}

\affil*[2]{\orgdiv{Institute of Mathematics and Informatics}, \orgname{Bulgarian Academy of Sciences}, \orgaddress{\street{8 Acad. G. Bontchev St.}, \city{Sofia} \postcode{1113}, \country{Bulgaria}}}

\abstract{The sum of independent, but not necessary identically distributed, exponential random variables follows a hypoexponential distribution. We focus on a particular case when all but one rate parameters of the
exponential variables are identical. This is known as exponentially modified Erlang distribution in molecular biology. We prove a
characterization of the exponential distribution, which complements previous characterizations via hypoexponential distribution with all rates different from each other.}



\maketitle

\section{Introduction and main results} 
Sums of exponentially distributed random variables play a central role in many stochastic models of real-world phenomena.  The {\it hypoexponential distribution} arises as a convolution of $n$ independent exponential distributions each with their own rate $\lambda _{i}$, the rate of the $i^{th}$ exponential distribution. It belongs to the class of {\it phase-type distributions}. Many processes
can be divided into sequential phases. If 
the time periods spent in different phases of the process are independent but not necessary identically distributed exponential variables, then  the overall time is hypoexponential. For example, the absorption time for a finite-state Markov chain follows this distribution.

\begin{figure}[h]
\centering
\includegraphics[width=0.5\textwidth]{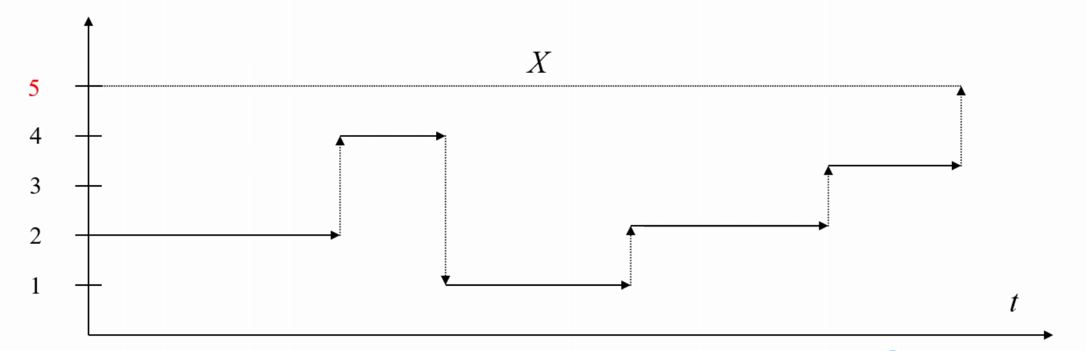}
\caption{Absorption time for a Markov chain with five states.}
\end{figure}

We will write $X_i\sim {\rm Exp}(\lambda_i)$, $\lambda_i>0$, if~$X_i$ for $i=1,2,\ldots$ has density
\be \label{exp_den}
f_i(x)=\lambda_i {\rm e}^{-\lambda_i x}, \quad x\ge 0 \quad \mbox{\it (exponential distribution)}.
\ee
The distribution of the sum 
$$ 
Y_n=X_1+X_2+\ldots +X_n \sim {\rm HypoE}(\lambda_1, \lambda_2, \ldots, \lambda_n),
$$
where $\lambda_i$ for $i=1,\,\ldots,\,n$ are not all identical, is called (general) {\it hypoexponential distribution} (e.g., \cite{LL19} and \cite{SKK16}). Assume that all $\lambda_i$'s are distinct, i.e., $\lambda _i\ne \lambda _j$ when $i\ne j$. It is well-known that under this condition, the density of 
$Y_n$ is given by (see \cite{R19}, p.309 and \cite{F71}, p.40, Problem 12) for $x>0$
\[
f_{Y_n}(x)=\sum_{j=1}^n \ \ell_j f_j(x),
         \qquad \lambda _i\ne \lambda _j, \  i\ne j.
\]
 Here the weight $\ell_j$ is defined as
$
\ell_j=\prod_{i=1, i\ne j}^n \lambda_i(\lambda_i-\lambda_j)^{-1}$. Thus, the density of the sum of independent exponential random variables with distinct parameters is linear combination of the individual densities.
For example, the density of $Y_2$
is 
\[
f_{Y_2}(x)
=\frac{\lambda_2}{\lambda_2-\lambda_1}f_1(x)+\frac{\lambda_1}{\lambda_1-\lambda_2}f_2(x), \quad \lambda_1\ne \lambda_2.
\]

\noindent Let $X_1$ and $X_2$ be two independent copies of a non-negative random variable 
$X$ and ${\mathsf E}[X]<\infty$. If $X\sim {\rm Exp}(\lambda)$, then 
$X_1+X_2/2 \sim {\rm HypoE}\left(\lambda, 2\lambda\right) $. It was proved in \cite{AV13} that this property of the exponential distribution is not shared by any other continuous distribution, i.e., 
 for $\lambda>0$
\be \label{AV}
  X\sim {\rm Exp}(\lambda) \quad \mbox{iff} \quad X_1+\frac{1}{2}X_2 \sim {\rm HypoE}\left(\lambda, 2\lambda\right).
\ee
The key argument in the proof is that the Laplace transform (LT) of the exponential distribution with $\lambda>0$
\be \label{solution}
\Phi(t)= \frac{\lambda}{\lambda + t}, \qquad t\ge 0
\ee
is the unique LT solution of the equation
\[
\Phi(t)\Phi\left(\frac{t}{2}\right)=2\Phi(t)-\Phi\left(\frac{t}{2}\right), \qquad t\ge 0.
\]
 Motivated by (\ref{AV}), in \cite{Y20} we extended it in two directions: 
 (i) for any number $n\ge 2$ of independent copies $X_1, X_2, \ldots, X_n$  of $X$, and 
 (ii) for the linear combination 
 $$
 \mu_1X_1+ \mu_2X_2+\ldots+\mu_nX_{n}\qquad \mu_i\ne \mu_j, \ i\ne j, \ \mu_i>0,
 $$
 i.e., with arbitrary positive and distinct coefficients $\mu_1, \mu_2,\,\ldots,\,\mu_n$. 
 Namely, it was proved in \cite{Y20}, under some additional assumptions,  that for  $\lambda>0$
\be \label{y20}
X\sim {\rm Exp}(\lambda) \quad \mbox{iff} \quad \sum_{k=1}^n\mu_kX_k
\sim {\rm HypoE}\left(\frac{\lambda}{\mu_1},\frac{\lambda}{\mu_2},\ldots, \frac{\lambda}{\mu_n}\right) .
\ee
This characterization was obtained by showing that (\ref{solution}) is the unique LT solution of the equation
\[ 
\Phi(\mu_1 t)\Phi(\mu_2 t)\cdots \Phi (\mu_n t) 
=\sum_{j=1}^n \bar{\ell}_j\Phi(\mu_j t), \qquad t\ge 0,
    \]
where $\bar{\ell}_j =  \prod_{i=1, i\ne j}^n \mu_j(\mu_j-\mu_i)^{-1}$. Thus, the case of the rate parameters $\lambda_i$'s in (\ref{exp_den}) being all different from each other was settled down. Note that characterization results in the case of distinct but not necessary  positive $\mu_1, \mu_2,\,\ldots,\,\mu_n$ were recently obtained, in \cite{RR23}, under an additional assumption.

 The other extreme case of all $\lambda_i$'s equal leads to Erlang distribution of the sum.
Assume $X_i\sim {\rm Exp}(\lambda)$, i.e., $\lambda_1=\lambda_2=\dots =\lambda_n=\lambda$ and let
\[
Y_n=X_1+X_2+\ldots + X_n.
\]
If $\Phi$ is the common LT of $X_i$, then for $t>0$
\be \label{Erl_eqn}
\Phi_{Y_n}(t)=\Phi^n(t) = \left(\frac{\lambda}{\lambda+t}\right)^n.
\ee
If we go in the opposite direction, assuming $Y_n\sim {\rm Erl}(n,\lambda)$, then (\ref{Erl_eqn}) yields
$\Phi_i(t)= \lambda(\lambda +t)^{-1}$ for each $i=1,2,\,\ldots,\,n$, which in turn implies $X_i\sim {\rm Exp}(\lambda)$. 
By words, if $X_i$ are independent and identically distributed r.v.'s 
and $Y_n$ is Erlang, then the common distribution is exponential.

 {\it The question arises whether a similar characterization holds when the rate parameters $\lambda_i$'s of ${\rm HypoE}(\lambda_1, \lambda_2, \ldots, \lambda_n)$ are neither all different nor all equal?}  It is our goal in this paper to show that, at least in one particular case,  the answer to this question is affirmative.

 Without the condition that all parameters $\lambda_i$'s are different or equal, the hypoexponential density has a quite complex form (\cite{JK03}). This makes the analysis of the general case difficult.  In this paper, we consider the particular case of "all-but-one-equal"
 rate parameters.
 More precisely, let $X_1, X_2,\ldots, X_{n+1}$ be independent copies of $X\sim {\rm Exp}(\lambda)$. 
Consider the sum 
\be \label{sum1}
X_1+ X_2+\ldots + X_n+wX_{n+1}, \qquad w>0,  w\ne 1.
\ee
This sum has a convoluted Erlang distribution, which is also known as  {\it exponentially modified Erlang (EME) distribution} (\cite{G16}, \cite{YFM17}).

Recall the well-known (e.g., \cite{B13}, p.240) Cram\'{e}r's condition. We say that
$X$ satisfies Cram\'{e}r's condition if there is a number $t_0>0$ such that ${\mathsf E}[{\rm e}^{tX}]<\infty$ for all $t\in (-t_0,t_0)$.
The next theorem establishes, under Cram\'{e}r's condition, a necessary and sufficient condition for $X\sim {\rm Exp}(\lambda)$.

{\bf Theorem} {\it Suppose that $ X_1, X_2,\,\ldots,\,X_{n+1}$, $n\ge 1$, are independent copies of a non-negative and absolutely continuous random variable $X$. Assume further that $X$ satisfies Cram\'{e}r's condition. Then for some $\lambda>0$, fixed positive integer $n$ and fixed positive real $w\ne 1$}
\be \label{thm2}
  X\sim {\rm Exp}(\lambda)\quad \mbox{iff} \quad \sum_{k=1}^{n} X_k + wX_{n+1}
\sim {\rm HypoE}\left(\lambda,\lambda, \ldots ,\lambda, \frac{\lambda}{w}\right) .
\ee

The hypoexponential family of distributions has found use in diverse applied fields, including queuing theory (\cite{BK11}), population genetics (\cite{SP01}), reliability analysis (\cite{KKK15}), medicine (\cite{DCF22}), and cell biology (\cite{YFM17}, \cite{GF19}). 
We focus on a particular member of this family, namely the exponentially modified Erlang distribution.  The (proper) Erlang distribution is applied in modeling the cell cycle phase progression as a series of sub-phase transitions with the same rate $\lambda$.   The relevant biological interpretation of the Erlang  model  is  that  each  cell  cycle phase can be viewed as a multi-step biochemical process that needs to be completed sequentially in order to advance to the next cell cycle phase.

Although the identical-stage model, 
is convenient from a mathematical perspective, it has been shown to be outperformed by a number of other distributions. Particularly, it was shown in \cite{G16} that one of the most appropriate distributions for representing cell cycle times is the EME distribution, which models a series of exponentially distributed random variables when  one of them has a different rate. Under this assumption the multi-stage cell cycle model is described as follows
\be \label{EME_system}
X_1 \overset{\lambda_1}{\longrightarrow} X_2 \overset{\lambda_1}{\longrightarrow}\cdots X_k\overset{\lambda_1}{\longrightarrow}X_{k+1}\overset{\lambda_2}{\longrightarrow}2X_1.
\ee
Note that, in the system (\ref{EME_system}), the rate of progression is identical through each of the initial $k$ stages of cell cycle and that we have added an additional exponentially distributed stage at the end whose rate, $\lambda_2$ is  distinct from the rate, $\lambda_1$, of the previous $k$ stages. Assume that the characteristic time (CT) of one of the transitions is significantly bigger than the CTs of all other transitions. Assume also that the CTs are independent. 
Then the overall cell cycle time  distribution may be approximated by the EME distribution capturing both the relatively rapid transitions and the slow transition (see \cite{G16}). Using the density convolution formula, for the density of (\ref{sum1}), it is not difficult to obtain
\[
f_{EME}(x)=\frac{\lambda}{w}e^{-\lambda x/w}\left(\frac{w}{w-1}\right)^n\left[ 1-\frac{\Gamma(n, (w-1)\lambda x/w}{(n-1)!}\right],
\]
where $\Gamma(n,t)=\int_t^\infty u^{n-1}e^{-u}\, du$ is the complementary incomplete gamma function.

In Section 2 we present some auxiliary results. We prove the  Theorem in Section~3. The last section includes some concluding remarks.

\section{Auxiliary results} 

\noindent Due to the independence assumption, the LT of (\ref{sum1}) equals $\Phi(wt)\Phi^n(t)$. If $\Phi$ is given by  (\ref{solution}), then $\Phi(wt)\Phi^n(t)$ is a product of linear
fractions and we can decompose it into sum of the Laplace transforms of $wX$ and $X$. 
Denote
\be \label{phi_1}
\Phi_1(t):=(w-1)\Phi(wt)\quad \mbox{and}\quad \Phi_2(t):=\frac{w-1}{w}\Phi(t).
\ee

{\bf Lemma 1}
The following identity holds
\beq \label{eqn_1_n}
\Phi_1(t)\Phi_2^n(t)
& = & \Phi_1(t)-\sum_{k=1}^{n} \Phi_2^{k}(t).
\eeq

{\bf Proof.} Without loss of generality assume that $X\sim {\rm Exp}(1)$. Note that  the following linear fraction decomposition holds
\be \label{decomp}
  \frac{w-1}{(1+wt)(1+t)}   =\frac{w}{1+wt}-\frac{1}{1+t}.
\ee
Recalling that $\Phi(t)=(1+t)^{-1}$, notation (\ref{phi_1}), and multiplying both sides of (\ref{decomp}) by $(w-1)/w$, we   obtain
\[
\Phi_1(t)\Phi_2(t) 
     = \Phi_1(t)-\Phi_2(t). 
\]
This proves (\ref{eqn_1_n}) for $n=1$.
Assuming that (\ref{eqn_1_n}) holds for $n$, for the $(n+1)^{th}$ term we have that
\nbeq
\Phi_1(t)\Phi_2^n(t)\Phi_2(t) 
    & = & 
    \left[\Phi_1(t)-\sum_{k=1}^{n} \Phi_2^{k}(t)\right]\Phi_2(t)
    \\
    & = & \Phi_1(t)\Phi_2(t)-\sum_{k=1}^{n} \Phi_2^{k+1}(t)
    \\
     & = &
\Phi_1(t)-\Phi_2(t)   
-\sum_{k=2}^{n+1} \Phi_2^{k}(t)
    \\ 
     & = &
     \Phi_1(t)
-\sum_{k=1}^{n+1} \Phi_2^{k}(t),
\neeq 
which proves (\ref{eqn_1_n}) for any $n$.

{\bf Lemma 2}
Let $n$ be any positive integer and $v\ne 1$ be a real number.

(i) For any integer $j\ge 1$ 
\be \label{le11}
\hspace{-0.3cm}
v\sum_{k=0}^{n-1} {k \choose j-1} v^{k} + (v-1)\sum_{k=0}^{n-1}{k \choose j} v^{k}={n \choose j}v^n.
\ee

(ii) For any integer $j\ge 2$ 
\be \label{lemma2}
\hspace{-1cm}\left(\frac{v}{v-1}\right)^{j-1}v\sum_{k=0}^{n-1} v^k +(v-1)\sum_{k=0}^{n-1}kv^k\ne nv^n . 
\ee   
{\bf Proof}. 
(i) The left-hand side of (\ref{le11}) is equivalent to 
\nbeq \label{le1}
\lefteqn{\sum_{k=0}^{n-1}{k \choose j-1} v^{k+1} + (v-1)\sum_{k=0}^{n-1} {k \choose j} v^{k} }\\
&  & = \sum_{k=0}^{n-1} \left[{k \choose j-1}+{k \choose j}\right]  v^{k+1}- \sum_{k=0}^{n-1}  {k \choose j}  v^{k}\nonumber \\
& & = \sum_{k=0}^{n-1} {k+1 \choose j} v^{k+1} - \sum_{k=0}^{n-1}  {k \choose j} v^{k} \nonumber \\
& & = \sum_{k=1}^{n} {k\choose j}  v^{k}- \sum_{k=1}^{n-1} {k \choose j} v^{k} \nonumber \\
&  & =  {n \choose j}v^n. \nonumber
\neeq

(ii) Using (\ref{le11}) with $j=1$, for the left-hand side of (\ref{lemma2}), we obtain
\nbeq
\lefteqn{\hspace{-1.5cm}\left(\frac{v}{v-1}\right)^{j-1}v\sum_{k=0}^{n-1} v^k -v\sum_{k=0}^{n-1} v^k+v\sum_{k=0}^{n-1} v^k+(v-1)\sum_{k=0}^{n-1}kv^k} \\
     & & = \left[\left(\frac{v}{v-1}\right)^{j-1}-1\right]v\sum_{k=0}^{n-1} v^k +nv^n \\
    & & = \left[ \left(\frac{v}{v-1}\right)^j-\frac{v}{v-1} \right](v^n-1)+ nv^n \ne nv^n.
\neeq
{\bf Remark.} It is not difficult to see that (\ref{le11}) can be generalized  to
\be \label{le11m}
v\sum_{k=0}^{n-1} {k+m \choose j-1} v^{k} + (v-1)\sum_{k=0}^{n-1}{k +m \choose j} v^{k}={n +m \choose j}v^n,
\ee
where $m$ is any non-negative integer.

\section{Proof of the theorem}

It follows from Lemma 1 that if  $X\sim {\rm Exp}(\lambda)$, then (\ref{eqn_1_n}) holds true. We will proceed with the proof of the opposite direction in the claim.
The case where $n=1$ is a particular case of (\ref{y20}) included in \cite{Y20}. Let $n\ge 2$. Consider the function $\Psi$ with the following series expansion 
\be \label{notation4}
\Psi(t):=\frac{1}{\Phi(t)}=\sum_{j=0}^\infty a_jt^j, \qquad t>0.
\ee
Note that, as a consequence of Cram\'{e}r's condition, the above series is uniformly convergent in a proper 
neighborhood of $t=0$ (see \cite{B13}, p.240).  To prove the theorem, 
it is sufficient to show that for some $\lambda>0$
\[
\Psi(t)=1+\lambda^{-1}t,
\]
i.e., the
coefficients of the series in (\ref{notation4}) are
\be \label{coefs}
a_0=1, \quad a_1=\lambda^{-1}>0, \quad a_j=0, \quad  j\ge 2.
\ee
Clearly,
\be \label{a0}
a_0=\Psi(0)=1. 
\ee
It follows from (\ref{phi_1}) and (\ref{eqn_1_n}) that 
\be \label{eqn_1_n2}
\frac{(w-1)^{n+1}}{w^n}\Phi(wt)\Phi^n(t)
=(w-1)\Phi(wt)-\sum_{k=1}^{n} \left(\frac{w-1}{w}\right)^{k}\Phi^{k}(t).
\ee
Dividing both sides of (\ref{eqn_1_n2}) by its left-hand side and changing the summation index, we obtain
\be \label{new_theorem21}
1=v^n \Psi^n(t)-(v-1)\Psi(wt)\sum_{k=0}^{n-1} v^k\Psi^k(t),
\ee
where, for notational simplicity, we set $v=w/(w-1)$.
To calculate the coefficients $a_j$ for $j\ge 1$, we 
differentiate both sides of (\ref{new_theorem21}) with respect to $t$ at $t=0$. After differentiating once at $t=0$ we have
\[
 \left[ nv^n -
v\sum_{k=0}^{n-1}v^k -(v-1)\sum_{k=0}^{n-1}kv^k\right]a_1 =0.
\]
It follows from (\ref{le11}) with $j=1$, that the coefficient in front of $a_1$ equals zero  and thus there exists a $\lambda>0$ such that 
\be \label{a1}
a_1=\lambda^{-1}.
\ee
Differentiating  (\ref{new_theorem21}) twice with respect to $t$ at $t=0$, we have
\nbeq
  \lefteqn{ \hspace{-1.5cm} \left[{n \choose 2}v^n -v\sum_{k=0}^{n-1}
            {k \choose 1}v^k-(v-1)\sum_{k=0}^{n-1} 
            {k \choose 2}v^k\right] a_1^2 }\\
            & & 
            + \left[nv^n- \left(\frac{v}{v-1}\right)v\sum_{k=0}^{n-1} v^k+ (v-1)\sum_{k=0}^{n-1}kv^k\right] a_2 =0.
\neeq
Lemma 2 with $j=2$ yields that the coefficient in front of $a_1^2$ is zero and the coefficient in front of $a_2$ is not zero. Therefore,
\be \label{a2}
a_2=0.
\ee
It remains to prove that $a_j=0$ for all $j\ge 3$. 
We will need the general Leibniz rule for differentiating {a} product of functions.
Denote by $y^{(j)}(x)$ the $j^{th}$ derivative of $y(x)$; $y^{(0)}(x):=y(x)$. 
Define a multi-index set
$\alphab=(\alpha_1,\alpha_2,\,\ldots,\,\alpha_{n})$ as a $n$-tuple of non-negative integers. Denote $ \| \alphab \|=\alpha_1+\alpha_2+\ldots+\alpha_n$ and $\Lambda_{j}:=\{\alphab : \|\alphab \|=j\}$.
The $j^{th}$ derivative (when exists) of the product $y_ 1(t)y_2(t)\cdots y_n(t)$ is given by (e.g. \cite{TL03})
\be \label{Lm}
\frac{{\rm d}^j}{{\rm d}t^j}\prod_{i=1}^n y_i(t)= \sum_{\Lambda_j} \left(\frac{j!}{\alpha_1!\alpha_2!\cdots \alpha_n!} \prod_{i=1}^n y_i^{(\alpha_i)}(t)\right). 
\ee
Let us write $\Lambda_{j}$ as union of three disjoint subsets as follows:
\[
\Lambda_{j}=\Lambda'_{j}\cup \Lambda''_{j}\cup \Lambda'''_{j},
\]
where 
\nbeq
\Lambda'_{j}  & =  &\{ \| \alphab \|=j : \ \mbox{only one of $\{\alpha_1, \alpha_2, \ldots, \alpha_n\}$ equals $j$ (others are zeros)} \}\\
\Lambda''_{j} &  =  & \{\|\alphab\|=j : \ \mbox{exactly $j$ of $\{\alpha_1, \alpha_2, \ldots, \alpha_n\}$ equal $1$ (others are zeros)}\} \\
\Lambda'''_{j} & =  &\{\|\alphab\|=j :\  \mbox{{there is an index $\alpha_i$ with $2\le \alpha_i<j$}} \}.
\neeq
Notice that by definition,  $\Lambda_j''$ is not empty only if  $j\le n$.

We will proceed by induction with respect to the index $j\ge 2$ of $a_j$. For $j=2$ we have already proved that $a_2=0$. Assuming $a_i=0$ for $2\le i\le j-1$, we will show that $a_j=0$.
 Since $a_0=1$, applying (\ref{Lm}), we obtain
\nbeq 
\frac{1}{j!}\frac{{\rm d}^j}{{\rm d}t^j} \Psi^n(t)\vert_{t=0} &  = & \sum_{\Lambda_{j}} \left(  \prod_{i=1}^n a_{\alpha_i} \right) = \sum_{\Lambda_j'} (\cdot) + \sum_{\Lambda_j''} (\cdot) + \sum_{\Lambda_j'''} (\cdot) \\
  &  =  &
    {n \choose 1}a_j a_0^{n-1} + {n \choose j}a_1^j a_0^{n-j} \\
    & = & na_j + {n \choose j}a_1^j.
\neeq
Notice that $\sum_{\Lambda_j'''} (\cdot)=0$ by the induction assumption.
Also, 
\nbeq
\lefteqn{\frac{1}{j!}\frac{{\rm d}^j}{{\rm d}t^j} \Psi(wt)\Psi^k(t)\vert_{t=0}  =  \sum_{\Lambda_{j}} \left( w^{\alpha_{k+1}} a_{\alpha_{k+1}} \prod_{i=1}^ka_{\alpha_i} \right)=\sum_{\Lambda'} (\cdot) + \sum_{\Lambda''} (\cdot) + \sum_{\Lambda'''} (\cdot) }\\
   & = &
   \left[ w^ja_0^ka_j + {k \choose 1}a_ja_0^{k-1}a_0\right] + \!  \left[{k \choose j-1}a_1^{j-1}a_0^{k-j+1}wa_1 + {k \choose j}a_1^j a_0^{k-j}a_0\right] \\
    & = &  \left[{k \choose j-1} w + {k \choose j}\right] a_1^j+\left( w^j+ k\right) a_j.
\neeq

\noindent Therefore, differentiating (\ref{new_theorem21}) $j$ times at $t=0$ and grouping the coefficients in front of $a_1^j$ and $a_j$, we write
\nbeq
  \lefteqn{ \hspace{-1.5cm} \left[{n \choose j}v^n - v\sum_{k=0}^{n-1} {k \choose j-1}v^k + 
    (v-1)\sum_{k=0}^{n-1}{k \choose j} v^k\right] a_1^j }\\
   & &  + 
    \left[nv^n-\left(\frac{v}{v-1}\right)^{j-1}v\sum_{k=0}^{n-1} v^k -(v-1)\sum_{k=0}^{n-1}kv^k\right] a_j =0.
\neeq
It then follows from Lemma 2 that the coefficient in front of $a_1^j$ is zero and the coefficient in front of $a_j$ is not zero. Therefore, for all $j\ge 2$
\be \label{aj}
a_j=0.
\ee
Now, (\ref{a0}), (\ref{a1}), (\ref{a2}), and (\ref{aj}) lead to  (\ref{coefs}), which completes the proof.

\section{Concluding remarks}
\label{sec:4}

In this paper we continue the study of the relation between the exponential and hypoexponential distributions, initiated in \cite{AV13} and extended in \cite{Y20}. The obtained characterization complements those in the above papers.   
Here we deal with a situation where the rate parameters $\lambda_i$'s in a convolution of exponential variables are not all different from each other.  The obtained result is of interest itself, however it can also serve as a basis for further investigations of more complex compositions of the rate parameters.

The exponential distribution is well known as a lifetime model in reliability theory. In particular, its use as a probability model for failure times of system's components is well justified. Thus, it is important to assess goodness-of-fit of the exponential distribution for a data set prior to applying the exponential model. Characterization results often serve as a useful device in obtaining goodness-of-fit tests. The presented here characterization can be used for testing the validity of a model based on the exponential distribution.



\bmhead{Acknowledgments}

The author acknowledges the valuable suggestions from the referees and the editor.

\bmhead{Statement of Conflict of Interest} The author states that there is no conflict of interest.

\end{document}